\documentclass{article}
\usepackage{amsfonts}

\def\G{{\cal{G}}}
\def\P{{\cal{P}}}

\newtheorem{lem}{Lemma}
\newtheorem{prp}{Proposition}
\newtheorem{cor}{Corollary}
\newtheorem{thm}{Theorem}
\newcommand{\qed  }{\hfill$\Box$\medskip}
\newcommand{\0}{\mbox{\boldmath $0$}}
\newcommand{\1}{\mbox{\boldmath $1$}}

\title
{\bf The graphs with all but two eigenvalues equal~to~$\pm1$}

\author
{Sebastian M. Cioab\u{a}
\\
{\it \small Department of Mathematical Sciences,}
\\
{\it \small University of Delaware, USA}
\\[3pt]
Willem H. Haemers\thanks{corresponding author; e-mail haemers@uvt.nl}
\\
{\it\small Department of Econometrics and Operations Research,}
\\
{\it\small Tilburg University, The Netherlands}
\\[3pt]
Jason Vermette
\\
{\it \small Department of Mathematical Sciences,}
\\
{\it \small University of Delaware, USA}
\\[3pt]
Wiseley Wong
\\
{\it \small Department of Mathematics,}
\\
{\it \small University of California, San Diego, USA}
}
\date{}

\begin{document}
\maketitle
\begin{abstract}\noindent
We determine all graphs whose adjacency matrix has
at most two eigenvalues (multiplicities included) different from $\pm 1$
and decide which of these graphs are determined by their spectrum.
This includes the so-called friendship graphs, which consist of a number
of edge-disjoint triangles meeting in one vertex.
It turns out that the friendship graph is determined by its spectrum,
except when the number of triangles equals sixteen.
\\[5pt]
\end{abstract}

\section{Introduction}

The friendship graph $F_k$ (also called Dutch windmill graph, or $k$-fan) consists of $k$ edge-disjoints triangles that meet in one vertex.
The famous friendship theorem (see Erd\H{o}s, R\'{e}nyi and S\'{o}s \cite{ERS} and Wilf \cite{W}) states that these are the only graphs
with the property that every pair of vertices contains a unique common neighbor (neighbors are called friends in the friendship theorem).
Clearly $F_k$ has $2k+1$ vertices and $3k$ edges, and $F_1=K_3$.
For convenience we shall assume that $k\geq 2$.
The adjacency matrix $A_k$ of $F_k$ has spectrum $\{\frac{1}{2}\pm\frac{1}{2}\sqrt{1+8k},\, 1^{k-1}, \, -1^k\}$
(multiplicities are denoted as exponents).
Wang, Belardo, Huang and Borovicanin~\cite{WBHB} conjectured that $F_k$ is determined by
the spectrum of the adjacency matrix $A_k$.
This conjecture caused some activity on the spectral characterization of $F_k$.
Das~\cite{D} claims to have a proof, but Abdollahi, Janbaz and Oboudi~\cite{AJO} found a mistake.
In addition these authors give correct proofs in some special cases.
In this paper, we prove that the conjecture from \cite{WBHB} is true if $k\neq 16$, and show that there is just one
counter example if $k=16$.

Although it has been conjectured by the second author that almost all graphs are
determined by the spectrum of the adjacency matrix, it is very often difficult to prove
the spectral characterization of a given graph (or family of graphs).
The spectrum of $A_k$ has two interesting properties that give much
information on the structure of the graph and bring a possible proof of the mentioned conjecture within reach.
The first property is that the second largest eigenvalue equals~$1$, and the second smallest eigenvalue is equal to~$-1$.
By eigenvalue interlacing (see for example~\cite{BH}, Section 2.5) it follows that every induced subgraph of a graph cospectral
with $F_k$ must have the second largest eigenvalue at most~1, and the second smallest eigenvalue at least~$-1$.
This gives a considerable reduction on the possible induced subgraphs (see Lemma~\ref{tabu-l}).
The second property is that $A_k^2-I$ has rank $2$ and is positive semi-definite.
This leads to conditions for the structure of $A_k^2$ (see Lemma~\ref{tools}).
Because of these observations we take a more general approach, and consider all graphs with the mentioned two properties.
Thus, we determine the graphs with only two eigenvalues $r$ and $s$ ($r>1$, $s<-1$) different from $\pm 1$.
We shall see that the disconnected ones have all components but one equal to $K_2$.
The connected ones come in three infinite families (one of which contains the friendship graphs)
and seven sporadic graphs.
No two non-isomorphic connected ones have the same spectrum,
but a disconnected graph can be cospectral and non-isomorphic to another one.
In particular, one of the sporadic graphs extended with some isolated edges is cospectral with $F_{16}$.

\section{Basics and tools}

We start with a well known result on equitable partitions (see for example~\cite{BH}).
Consider a partition ${\P}=\{V_1,\ldots,V_m\}$ of the set $V=\{1,\ldots,n\}$.
The characteristic matrix $\chi_\P$ of $\P$ is the $n\times m$ matrix whose
columns are the character vectors of $V_1,\ldots,V_m$.
Consider a symmetric matrix $A$ of order $n$, with rows and columns partitioned
according to $\P$.
The partition of $A$ is {\em equitable} if each submatrix $A_{i,j}$ formed by the
rows of $V_i$ and the columns of $V_j$ has constant row sums $q_{i,j}$.
The $m\times m$ matrix $Q=(q_{i,j})$ is called the {\em quotient matrix}
of $A$ with respect to $\P$.
\begin{lem}\label{partition}
The matrix $A$ has the following two kinds of eigenvectors and eigenvalues:
\begin{itemize}
\item[(i)]
The eigenvectors in the column space of $\chi_\P$; the corresponding eigenvalues coincide with the eigenvalues of $Q$.
\item[(ii)]
The eigenvectors orthogonal to the columns of $\chi_\P$; the corresponding eigenvalues of $A$ remain unchanged if
some scalar multiple of the all-one block $J$ is added to block $A_{i,j}$ for each $i,j\in\{1,\ldots,m\}$.
\end{itemize}
\end{lem}
The reverse identity matrix of order $n$ is denoted by $R_n$.
Thus $R_{2k}$ is the adjacency matrix of $kK_2$, the disjoint union of $k$ edges.
We illustrate the use of Lemma~\ref{partition} with an example.
Consider the following partitioned matrix $A$ with quotient matrix $Q$:
\[ 
A =
\left[\begin{array}{cc}
J-I_{a}& J \\ J & R_{2k}
\end{array}\right]
\ ,\ Q =
\left[\begin{array}{cc}
a-1 & 2k \\ a & 1
\end{array}\right].
\]
The eigenvalues of $Q$ are $(a\pm\sqrt{a^2+8ak-4a+4})/2$, so they are also eigenvalues of $A$.
The other eigenvalues of $A$ remain the same if we subtract $J$ from the blocks equal to $J$ or $J-I_a$.
Then $A$ and $Q$ become
\[ 
A' =
\left[\begin{array}{cc}
-I_a & O \\ O & R_{2k}
\end{array}\right]
\ ,\ Q' =
\left[\begin{array}{cc}
-1 & 0 \\ 0 & 1
\end{array}\right].
\]
The part of the spectrum of $A'$, which is not in the spectrum of $Q'$ is $\{1^{k-1},\, -1^{k+a-1}\}$.
Thus we find that $A$ has spectrum $\{(a\pm\sqrt{a^2+8ak-4a+4})/2,\, 1^{k-1},\, -1^{k+a-1}\}$.
We see that a graph $G$ with adjacency matrix $A$ belongs to the classification.
Note that if $a=1$, then $G$ is the friendship graph $F_k$.
\begin{prp}\label{start}
Let $G$ be a graph with $n$ vertices and adjacency matrix $A$.
\begin{itemize}
\item[$(i)$]
If $A$ has all its eigenvalues equal to $\pm 1$, then $G=\frac{n}{2}K_2$.
\item[$(ii)$]
If $A$ has all but one eigenvalue equal to $\pm 1$, then $G$ is the disjoint union of
complete graphs with all but one connected components equal to $K_2$.
\item[$(iii)$]
If $A$ has just two eigenvalues, $r$ and $s$ ($r\geq s$) different from $\pm 1$,
then $r>1$ and $s<-1$, or $G$ is a disjoint union of complete graphs
with two connected components different from $K_2$.
\end{itemize}
\end{prp}
{\bf Proof}.
If $A$ has an eigenvalue $s<-1$, then the largest eigenvalue of $A$ is greater than $1$
(by the Perron-Frobenius theorem); this case is captured by the first option of $(iii)$.
If $A$ has smallest eigenvalue at least $-1$, then $G$ is the disjoint union of cliques
(see for example~\cite{vDH}), which leads to the other possibilities.
\qed
\begin{lem}\label{tools}
Suppose $r>1$ and $s<-1$ are the only eigenvalues of $G$ different from $\pm 1$.
\begin{itemize}
\item[$(i)$]
One connected component of $G$ has all vertex degrees at least $2$, and all other connected components are isomorphic to $K_2$.
\item[$(ii)$]
Suppose $u$ and $v$ are distinct vertices with degrees $d_u$ and $d_v$, and each neighbor of $u$ is also a neighbor of $v$.
Then $d_v-d_u \geq 3$.
\end{itemize}
\end{lem}
{\bf Proof}.
$(i)$
Suppose $u$ is a vertex of degree $1$.
Let $v$ be the neighbor of $u$, and assume that $v$ has another neighbor $w$ of degree $d_w$.
The $2\times 2$ principal submatrix of $A^2-I$ corresponding to $u$ and $w$ equals
\[ S = \left[\begin{array}{cc} 0 & 1\\1 & d_w-1\end{array}\right] .
\]
We have $\det S<0$, whilst $A^2-I$ is positive semi-definite.
This is a contradiction proving that $v$ has degree 1.
\\
$(ii)$
The $2\times 2$ principal submatrix of $A^2-I$ corresponding to $u$ and $v$ equals
\[ S = \left[\begin{array}{cc} d_u-1 & d_u\\d_u & d_v-1\end{array}\right] .
\]
If $d_v\leq d_u +2$, then $\det S \leq (d_u-1)(d_u+1)-d_u^2 < 0$, contradiction.
\qed
\\
Note that $(ii)$ of Lemma~\ref{tools} implies that two vertices $u$ and $v$ cannot have the same set of neighbors.

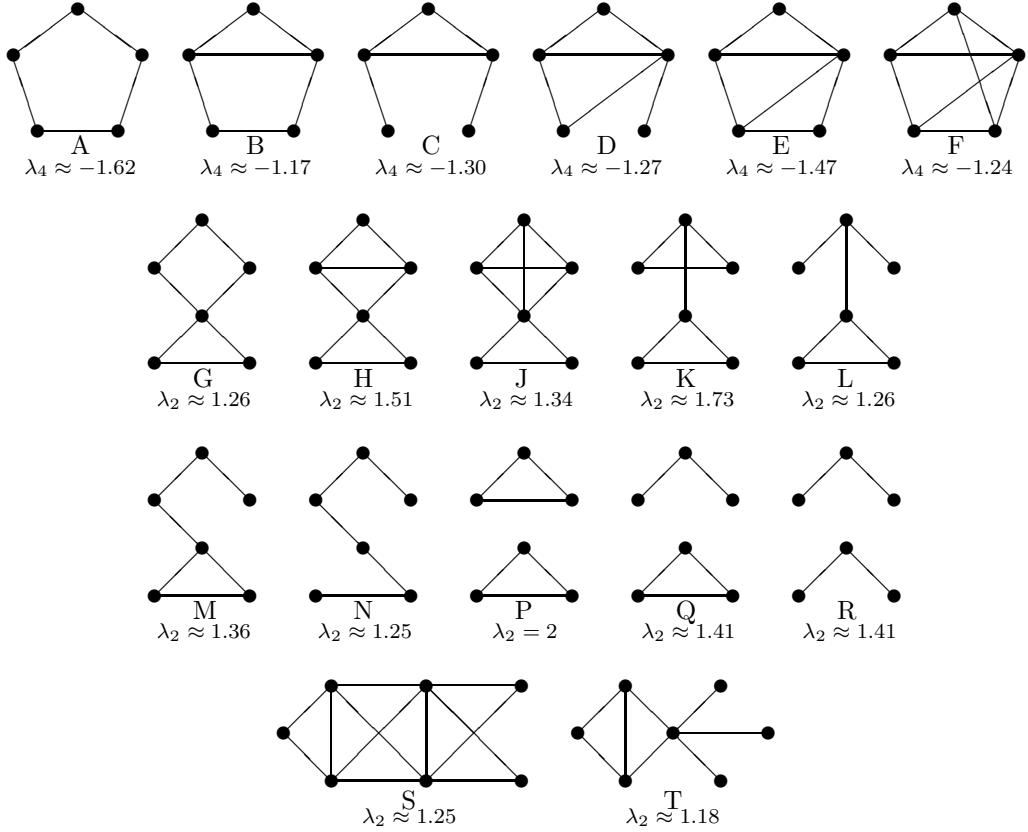
\begin{figure}[h]
\begin{center}
\setlength{\unitlength}{1.8pt}
\begin{picture}(37,36)(-2,0)
\put(5,9){\circle*{3}}
\put(22,9){\circle*{3}}
\put(0,25){\circle*{3}}
\put(27,25){\circle*{3}}
\put(13.5,34.5){\circle*{3}}
\put(5,9){\line(1,0){17}}
\put(5,9){\line(-1,3){5}}
\put(22,9){\line(1,3){5}}
\put(13.5,35){\line(4,-3){13.5}}
\put(13.5,35){\line(-4,-3){13.5}}
\put(12,4){A}
\put(2.25,0){\footnotesize $\lambda_4\approx-1.62$}
\end{picture}
\begin{picture}(35,36)(0,0)
\put(5,9){\circle*{3}}
\put(22,9){\circle*{3}}
\put(0,25){\circle*{3}}
\put(27,25){\circle*{3}}
\put(13.5,34.5){\circle*{3}}
\put(0,25){\line(1,0){27}}
\put(5,9){\line(1,0){17}}
\put(5,9){\line(-1,3){5}}
\put(22,9){\line(1,3){5}}
\put(13.5,35){\line(4,-3){13.5}}
\put(13.5,35){\line(-4,-3){13.5}}
\put(12,4){B}
\put(2.25,0){\footnotesize $\lambda_4\approx-1.17$}
\end{picture}
\begin{picture}(35,36)(0,0)
\put(5,9){\circle*{3}}
\put(22,9){\circle*{3}}
\put(0,25){\circle*{3}}
\put(27,25){\circle*{3}}
\put(13.5,34.5){\circle*{3}}
\put(0,25){\line(1,0){27}}
\put(5,9){\line(-1,3){5}}
\put(22,9){\line(1,3){5}}
\put(13.5,35){\line(4,-3){13.5}}
\put(13.5,35){\line(-4,-3){13.5}}
\put(12,4){C}
\put(2.25,0){\footnotesize $\lambda_4\approx-1.30$}
\end{picture}
\begin{picture}(35,36)(0,0)
\put(5,9){\circle*{3}}
\put(22,9){\circle*{3}}
\put(0,25){\circle*{3}}
\put(27,25){\circle*{3}}
\put(13.5,34.5){\circle*{3}}
\put(0,25){\line(1,0){27}}
\put(5,9){\line(4,3){22}}
\put(5,9){\line(-1,3){5}}
\put(22,9){\line(1,3){5}}
\put(13.5,35){\line(4,-3){13.5}}
\put(13.5,35){\line(-4,-3){13.5}}
\put(12,4){D}
\put(2.25,0){\footnotesize $\lambda_4\approx-1.27$}
\end{picture}
\begin{picture}(35,36)(0,0)
\put(5,9){\circle*{3}}
\put(22,9){\circle*{3}}
\put(0,25){\circle*{3}}
\put(27,25){\circle*{3}}
\put(13.5,34.5){\circle*{3}}
\put(0,25){\line(1,0){27}}
\put(5,9){\line(1,0){17}}
\put(5,9){\line(4,3){22}}
\put(5,9){\line(-1,3){5}}
\put(22,9){\line(1,3){5}}
\put(13.5,35){\line(4,-3){13.5}}
\put(13.5,35){\line(-4,-3){13.5}}
\put(12,4){E}
\put(2.25,0){\footnotesize $\lambda_4\approx-1.47$}
\end{picture}
\begin{picture}(35,36)(0,0)
\put(5,9){\circle*{3}}
\put(22,9){\circle*{3}}
\put(0,25){\circle*{3}}
\put(27,25){\circle*{3}}
\put(13.5,34.5){\circle*{3}}
\put(0,25){\line(1,0){27}}
\put(5,9){\line(1,0){17}}
\put(5,9){\line(4,3){22}}
\put(5,9){\line(-1,3){5}}
\put(22,9){\line(1,3){5}}
\put(13.5,35){\line(4,-3){13.5}}
\put(13.5,35){\line(-4,-3){13.5}}
\put(13.5,35){\line(1,-3){8.5}}
\put(12,4){F}
\put(2.25,0){\footnotesize $\lambda_4\approx-1.24$}
\end{picture}
\\[15pt]
\hspace{10pt}
\begin{picture}(32,40)(0,0)
\put(1,9){\circle*{3}}
\put(1,29){\circle*{3}}
\put(21,9){\circle*{3}}
\put(21,29){\circle*{3}}
\put(11,19){\circle*{3}}
\put(11,39){\circle*{3}}
\put(1,9){\line(1,1){20}}
\put(1,9){\line(1,0){20}}
\put(1,29){\line(1,1){10}}
\put(21,29){\line(-1,1){10}}
\put(21,9){\line(-1,1){20}}
\put(9,4){G}
\put(1.5,0){\footnotesize $\lambda_2\approx 1.26$}
\end{picture}
\begin{picture}(32,40)(0,0)
\put(1,9){\circle*{3}}
\put(1,29){\circle*{3}}
\put(21,9){\circle*{3}}
\put(21,29){\circle*{3}}
\put(11,19){\circle*{3}}
\put(11,39){\circle*{3}}
\put(1,9){\line(1,1){20}}
\put(1,9){\line(1,0){20}}
\put(1,29){\line(1,1){10}}
\put(1,29){\line(1,0){20}}
\put(21,29){\line(-1,1){10}}
\put(21,9){\line(-1,1){20}}
\put(9,4){H}
\put(1.75,0){\footnotesize $\lambda_2\approx 1.51$}
\end{picture}
\begin{picture}(32,40)(0,0)
\put(1,9){\circle*{3}}
\put(1,29){\circle*{3}}
\put(21,9){\circle*{3}}
\put(21,29){\circle*{3}}
\put(11,19){\circle*{3}}
\put(11,39){\circle*{3}}
\put(1,9){\line(1,1){20}}
\put(1,9){\line(1,0){20}}
\put(11,19){\line(0,1){20}}
\put(1,29){\line(1,1){10}}
\put(1,29){\line(1,0){20}}
\put(21,29){\line(-1,1){10}}
\put(21,9){\line(-1,1){20}}
\put(9,4){J}
\put(1.5,0){\footnotesize $\lambda_2\approx 1.34$}
\end{picture}
\begin{picture}(32,40)(0,0)
\put(1,9){\circle*{3}}
\put(1,29){\circle*{3}}
\put(21,9){\circle*{3}}
\put(21,29){\circle*{3}}
\put(11,19){\circle*{3}}
\put(11,39){\circle*{3}}
\put(1,9){\line(1,1){10}}
\put(1,9){\line(1,0){20}}
\put(11,19){\line(0,1){20}}
\put(1,29){\line(1,1){10}}
\put(1,29){\line(1,0){20}}
\put(21,29){\line(-1,1){10}}
\put(21,9){\line(-1,1){10}}
\put(9,4){K}
\put(1.5,0){\footnotesize $\lambda_2\approx 1.73$}
\end{picture}
\begin{picture}(32,40)(0,0)
\put(1,9){\circle*{3}}
\put(1,29){\circle*{3}}
\put(21,9){\circle*{3}}
\put(21,29){\circle*{3}}
\put(11,19){\circle*{3}}
\put(11,39){\circle*{3}}
\put(1,9){\line(1,1){10}}
\put(1,9){\line(1,0){20}}
\put(11,19){\line(0,1){20}}
\put(1,29){\line(1,1){10}}
\put(21,29){\line(-1,1){10}}
\put(21,9){\line(-1,1){10}}
\put(9,4){L}
\put(1.5,0){\footnotesize $\lambda_2\approx 1.26$}
\end{picture}
\\[15pt]
\hspace{10pt}
\begin{picture}(32,40)(0,0)
\put(1,9){\circle*{3}}
\put(1,29){\circle*{3}}
\put(21,9){\circle*{3}}
\put(21,29){\circle*{3}}
\put(11,19){\circle*{3}}
\put(11,39){\circle*{3}}
\put(1,9){\line(1,1){10}}
\put(1,9){\line(1,0){20}}
\put(1,29){\line(1,1){10}}
\put(21,29){\line(-1,1){10}}
\put(21,9){\line(-1,1){20}}
\put(9,4){M}
\put(1.5,0){\footnotesize $\lambda_2\approx 1.36$}
\end{picture}
\begin{picture}(32,40)(0,0)
\put(1,9){\circle*{3}}
\put(1,29){\circle*{3}}
\put(21,9){\circle*{3}}
\put(21,29){\circle*{3}}
\put(11,19){\circle*{3}}
\put(11,39){\circle*{3}}
\put(1,9){\line(1,0){20}}
\put(1,29){\line(1,1){10}}
\put(21,29){\line(-1,1){10}}
\put(21,9){\line(-1,1){20}}
\put(9,4){N}
\put(1.5,0){\footnotesize $\lambda_2\approx 1.25$}
\end{picture}
\begin{picture}(32,40)(0,0)
\put(1,9){\circle*{3}}
\put(1,29){\circle*{3}}
\put(21,9){\circle*{3}}
\put(21,29){\circle*{3}}
\put(11,19){\circle*{3}}
\put(11,39){\circle*{3}}
\put(1,9){\line(1,1){10}}
\put(1,9){\line(1,0){20}}
\put(1,29){\line(1,1){10}}
\put(1,29){\line(1,0){20}}
\put(21,29){\line(-1,1){10}}
\put(21,9){\line(-1,1){10}}
\put(9,4){P}
\put(4.25,0){\footnotesize $\lambda_2 = 2$}
\end{picture}
\begin{picture}(32,40)(0,0)
\put(1,9){\circle*{3}}
\put(1,29){\circle*{3}}
\put(21,9){\circle*{3}}
\put(21,29){\circle*{3}}
\put(11,19){\circle*{3}}
\put(11,39){\circle*{3}}
\put(1,9){\line(1,1){10}}
\put(1,9){\line(1,0){20}}
\put(1,29){\line(1,1){10}}
\put(21,29){\line(-1,1){10}}
\put(21,9){\line(-1,1){10}}
\put(9,4){Q}
\put(1.75,0){\footnotesize $\lambda_2\approx 1.41$}
\end{picture}
\begin{picture}(32,40)(0,0)
\put(1,9){\circle*{3}}
\put(1,29){\circle*{3}}
\put(21,9){\circle*{3}}
\put(21,29){\circle*{3}}
\put(11,19){\circle*{3}}
\put(11,39){\circle*{3}}
\put(1,9){\line(1,1){10}}
\put(1,29){\line(1,1){10}}
\put(21,29){\line(-1,1){10}}
\put(21,9){\line(-1,1){10}}
\put(9,4){R}
\put(1.75,0){\footnotesize $\lambda_2\approx 1.41$}
\end{picture}
\\[15pt]
\hspace{40pt}
\begin{picture}(60,30)(-7,0)
\put(0,9){\circle*{3}}
\put(20,9){\circle*{3}}
\put(0,29){\circle*{3}}
\put(20,29){\circle*{3}}
\put(40,9){\circle*{3}}
\put(40,29){\circle*{3}}
\put(-10,19){\circle*{3}}
\put(0,9){\line(1,1){20}}
\put(0,9){\line(1,0){40}}
\put(0,9){\line(0,1){20}}
\put(20,9){\line(1,1){20}}
\put(20,9){\line(-1,1){20}}
\put(20,9){\line(0,1){20}}
\put(20,29){\line(1,-1){10}}
\put(0,29){\line(1,0){40}}
\put(20,29){\line(1,-1){20}}
\put(-10,19){\line(1,-1){10}}
\put(-10,19){\line(1,1){10}}
\put(14.75,3){S}
\put(6.5,0){\footnotesize $\lambda_2\approx 1.25$}
\end{picture}
\begin{picture}(60,30)(-7,0)
\put(0,9){\circle*{3}}
\put(20,9){\circle*{3}}
\put(0,29){\circle*{3}}
\put(10,19){\circle*{3}}
\put(30,19){\circle*{3}}
\put(20,29){\circle*{3}}
\put(-10,19){\circle*{3}}
\put(0,9){\line(1,1){20}}
\put(0,9){\line(0,1){20}}
\put(0,29){\line(1,-1){20}}
\put(-10,19){\line(1,-1){10}}
\put(-10,19){\line(1,1){10}}
\put(10,19){\line(1,0){20}}
\put(8,3){T}
\put(0,0){\footnotesize $\lambda_2\approx 1.18$}
\end{picture}
\caption{Forbidden induced subgraphs}\label{tabu-f}
\end{center}
\end{figure}

Define $\G$ to be the set of connected graphs with eigenvalues $r>1$ and
$s<-1$, and all other eigenvalues equal to $\pm 1$.
By the above results, in order to find all graphs with at most two eigenvalues different from $\pm 1$,
it suffices to determine $\G$.
We start with a list of forbidden induced subgraphs.
\begin{lem}\label{tabu-l}
No graph in $\G$ has one of the graphs presented in Figure~\ref{tabu-f} is an induced subgraph.
\end{lem}
{\bf Proof}. Each graph in Figure~\ref{tabu-f} has its second largest eigenvalue $\lambda_2$ strictly greater than $1$,
or its second smallest eigenvalue $\lambda_{n\!-\!1}$ strictly less than $-1$. Interlacing completes the proof.~\qed
\section{Main results}
We begin with the description of the graphs in $\G$.
The proof will be given in the next section.
\begin{thm}\label{main}
The adjacency matrices and spectra of the graphs in $\G$ are as follows:
\begin{itemize}
\item[$(i)$]
$\left[\begin{array}{cc}
O &J-I_{m} \\ J-I_{m} & O
\end{array}\right]
\ (m\geq 3)$
\\[2pt]
with spectrum $ \{\pm(m-1),\ 1^{m-1},\ -1^{m-1} \} $,
\item[$(ii)$]
$\left[\begin{array}{cc}
J-I_a & J \\
J     & R_{2k}
\end{array}\right]
\ (a\geq 1,\ k\geq 2 )$
\\[2pt]
with spectrum $\{ \frac{a}{2}\pm\frac{1}{2}\sqrt{a^2+8ak-4a+4},\ 1^{k-1},\ -1^{a+k-1} \} $,
\item[$(iii)$]
$\left[\begin{array}{cc}
R_{2\ell} & J \\
J      & R_{2m}
\end{array}\right]
\ (\ell\geq m\geq 2)$
\\[2pt]
with spectrum $\{ 1 \pm 2\sqrt{\ell m},\ 1^{\ell+m-2},\ -1^{\ell+m} \}$,
\item[$(iv)$]
$\left[\begin{array}{cc}
O &N \\ N^\top & O
\end{array}\right]$
where
$N= \left[\begin{array}{cc}1 & \1^\top \\ \1 & I_4\end{array}\right]$,  or
$N= \left[\begin{array}{cc}J-I_3 & J \\ O & J-I_3\end{array}\right]$
\\[2pt]
with spectra $\{\pm 3,\ 1^4,\ -1^4 \}$, and $\{\pm 4,\ 1^5,\ -1^5 \}$, respectively,
\item[$(v)$]
$ \left[\begin{array}{ccc}
J-I_a   & J       & \1 \\
J       & J-I_{b} & \0 \\
\1^\top & \0^\top & 0
\end{array}\right]$ where $(a,b)=(6,5)$, $(4,6)$, or $(3,8)$
\\[2pt]
with spectra $\{4\pm2\sqrt{10},\ 1^1,\ -1^9 \}$, $\{(7\pm\sqrt{129})/2,\ 1^1,\ -1^{8}\}$,
$\{4\pm\sqrt{37},\ 1^1,\ -1^{9}\}$,
\item[$(vi)$]
$\left[
\begin{array}{ccc}
J-I_a & J     & O \\
J     & O     & J-I_m \\
O     & J-I_m & O
\end{array}\right]$ where $(a,m)=(3,5)$ or $(4,4)$
\\[2pt]
with spectra
$\{(1\pm\sqrt{129})/2,\ 1^{5},\ -1^{6} \}$, $\{1\pm2\sqrt{7},\ 1^4,\ -1^6 \}$.
\end{itemize}
\end{thm}
We see that $\G$ contains three infinite families and seven sporadic graphs.
From the given spectra it follows straightforwardly that
\begin{cor}
No two graphs in $\G$ are cospectral.
\end{cor}
\begin{thm}
Suppose $G$ and $G'$ are nonisomorphic cospectral graphs with at most two eigenvalues
different from $\pm 1$.
Then $G=H+\alpha K_2$ and $G'=H'+\alpha' K_2$, where $H$ and $H'$ are one of the following pairs of graphs in $\G$:
\begin{itemize}
\item
Both $H$ and $H'$ are of type $(iii)$ with parameters $(\ell,m)$
and $(\ell',m')$, where $\ell m=\ell' m'$, 
\item
One is of type $(iii)$ with $\ell,m\geq 2$, and the other of type $(ii)$ with $a=2$ and $k=\ell m$,
\item
One is of type $(iv)$ and the other one of type $(i)$ with $m=4$, or $5$, respectively,
\item
One is of type $(ii)$ with $(a,k)=(1,16)$ or $(2,7)$, and the other of type $(vi)$ with $(a,m)=(3,5)$, or $(4,4)$, respectively.
\end{itemize}
\end{thm}
{\bf Proof.}
The disjoint union of complete graphs in known to be determined by its spectrum (see~\cite{vDH}).
So, by Proposition~\ref{start} and Lemma~\ref{tools}$(i)$, $G$ and $G'$ must have the described form.
Next observe that $H$ and $H'$ share the eigenvalues $r>1$ and $s<-1$.
Using this we easily find the given possibilities for $H$ and $H'$.
\qed.
\\[3pt]
It we take $\alpha'\geq\alpha=0$ we find the graphs in $\G$ having a nonisomorphic cospectral mate.
\begin{cor}
A graph $G\in\G$ is determined by its spectrum, unless $G$ is one of the following
\begin{itemize}
\item
$G$ is of type $(ii)$  and $(a,k)=(1,16)$ or $(2,7)$,
\item
$G$ is of type $(ii)$ with $a=2$ and $k$ a composite number,
\item
$G$ is of type $(iii)$ and $\ell m$ has a divisor strictly between $\ell$ and $m$,
\item
$G$ is of type $(iv)$.
\end{itemize}
\end{cor}
Thus we have that the friendship graph $F_k$, which is Case $(ii)$ with $a=1$, is determined by its spectrum,
except when $k=16$.
The friendship graph $F_{16}$ is cospectral with $G + 10 K_2$, where $G\in\G$ is of type $(vi)$ with $(a,m)=(3,5)$.

\section{The proof}
Here we give the proof of Theorem~\ref{main}.

In all cases we see that the corresponding quotient matrix has two eigenvalues different from
$\pm 1$ and with Lemma~\ref{partition} it straightforwardly follows that the remaining eigenvalues of
the graph are all equal to $\pm 1$. So all graphs of Theorem~\ref{main} are in $\G$.

E.R. van Dam and E. Spence~\cite{vDS} classified all bipartite graphs with four distinct eigenvalues.
Their Proposition~8 gives the bipartite graphs in $\G$, described in $(i)$ and $(iv)$.

In the remainder of the proof it is assumed that $G\in\G$ is not bipartite.
We define $C$ to be a clique in $G$ with maximum size.
By Lemma~\ref{tabu-l} (graphs A and N) $G$ contains no induced odd cycles of length
five or more, therefore $|C|\geq 3$.
If there are more than one cliques of maximum size, we choose one for which the number of outgoing edges is minimal.
The following lemma is the key to our approach.
\begin{lem}
The vertex set of $C$ can be partitioned into two nonempty subsets $X$ and $Y$ (say), such
that the neighborhood of any vertex outside $C$ intersects $C$ in $X$, $Y$, or $\emptyset$.
\end{lem}
{\bf Proof}.
If $|C|=n-1$ the result is obvious.
So assume $3\leq |C|\leq n-2$.
Take vertices $x$ and $y$ outside $C$, and let $X$ and $Y$ consist of the neighbors of $x$ and $y$ in $C$, respectively.
Note that $X$ and $Y$ are proper subsets of $C$, since otherwise $C$ is not maximal.
Suppose that $X\cap Y\neq\emptyset$ but $X\not\subset Y$. 
Then there exist vertices $u\in X\cap Y$ and $v\in X\setminus Y$.
Let $w$ be a vertex in $C\setminus X$. Then the subgraph induced by $\{u,v,w,x,y\}$ is a forbidden subgraph D, E, or F.
Therefore, if $X$ and $Y$ are not disjoint, then $X\subset Y$, and analogously $Y\subset X$. 
Thus $X\cap Y\neq\emptyset$ implies $X=Y$.
%
If $X\cap Y=\emptyset$, assume there exist vertices $u\in X$, $v\in Y$, and $z\in C\setminus(X\cup Y)$,
then $\{z,u,v,x,y\}$ induces a forbidden subgraph B or C.
This implies that if $X$ and $Y$ are disjoint and both nonempty, then $X\cup Y = C$.
\qed\\
Let $\Gamma X$ and $\Gamma Y$ denote the set of vertices outside $C$ adjacent to $X$ and $Y$ respectively.
The set of vertices not adjacent to any vertex of $C$ will be denoted by $\Omega$.
Some of these sets may be empty, but clearly $\Gamma X$ or $\Gamma Y$ is nonempty
(otherwise $G$ would be disconnected or complete).
We choose $\Gamma X\neq\emptyset$ and distinguish three cases:
(1) both $\Gamma Y$ and $\Omega$ are empty, (2) only $\Omega$ is empty, and
(3) $\Omega$ is nonempty.
For convenience we define $a=|X|$, $b=|Y|$, and $c=|C|=a+b$.

\subsection{$\Gamma Y$ and $\Omega$ are empty}
Assume $b=1$.
Then $\Gamma X$ contains no edges, because $C$ is maximal.
The vertex $v\in C\setminus X$ and a vertex in $\Gamma X$ are nonadjacent with the same neighbors,
which is impossible by Lemma~\ref{tools}$(ii)$.
Therefore $b\geq 2$.
Choose two vertices $u$ and $v$ from $Y$, and choose $w\in X$.
Suppose $x\in\Gamma X$ has two neighbors $y$ and $z$ in $\Gamma X$, then
$\{u,v,w,x,y,z\}$ induces graph J or $\{v,w,x,y,z\}$ induces graph D from Figure~\ref{tabu-f},
therefore any vertex $x\in\Gamma X$ has at most one neighbor in $\Gamma X$.
By Lemma~\ref{tools}$(ii)$, it is not possible that $x\in\Gamma X$ has one neighbor in
$\Gamma X$ and $y\in\Gamma X$ has no neighbor in $\Gamma X$.
We conclude that either all vertices of $\Gamma X$ have exactly one neighbor in $\Gamma X$,
or $\Gamma X$ contains no edges.
In the first case $G$ has the following adjacency matrix $A$ with quotient matrix $Q$:
\[ A=\left[\begin{array}{ccc}
J-I_a & J & J \\
J & J-I_{b} & O \\
J & O & R_{n-c}
\end{array}\right],\ 
Q=\left[\begin{array}{ccc}
a-1 & b   & n-c \\
a   & b-1 & 0 \\
a   & 0   & 1
\end{array}\right].
\]
Computing $\det(Q+I)$ and $\det(Q-I)$ shows that $Q$ has no eigenvalue $-1$, and
$Q$ has an eigenvalue $1$ if and only if $b=2$.
In case $b=2$ we can rewrite $A$ as
\[ A=\left[\begin{array}{cc}
J-I_a & J \\
J     & R_{2k}
\end{array}\right]
\]
with $k\geq 2$.
Thus we obtained the graphs of Case~$(ii)$.

If $\Gamma X$ has no edges and at least two vertices, then these two vertices have the same neighbors,
contradiction.
So $|\Gamma X|=1$ and we find
\[
A=\left[
\begin{array}{ccc}
J-I_a   & J       & \1 \\
J       & J-I_{b} & \0 \\
\1^\top & \0^\top & 0
\end{array}\right],\ 
Q=\left[\begin{array}{ccc}
a-1 & b & 1 \\
a & b-1 & 0 \\
a & 0 & 0
\end{array}\right] .
\]
The quotient matrix $Q$ has no eigenvalue $-1$ and an eigenvalue $1$ if only if
$(a,b)=(6,5)$, $(4,6)$, or $(3,8)$, which leads to Case~$(v)$.

\subsection{$\Gamma X$ and $\Gamma Y$ are nonempty, and $\Omega$ is empty}

{\bf Claim~1~} $a\leq 2\ \mbox{ or } b\leq 2$.
\\[3pt]
{\bf Proof}.
Suppose $a\geq b\geq 3$ and suppose $\{x,y\}$ is an edge in $\Gamma Y$.
Let $u,v,w$ be three distinct vertices in $X$, and choose $z\in Y$.
Then $\{u,v,w,x,y,z\}$ induce graph~J from Lemma~\ref{tabu-l}.
So $\Gamma Y$ contains no edges.
Similarly $\Gamma X$ has no edges.
Now forbidden subgraph~S from Lemma~\ref{tabu-l} implies that
a vertex in $\Gamma X$ is adjacent to all, or all but one vertices in
$\Gamma Y$ (and vice versa).
Let $x$ be a vertex in $\Gamma X$ and suppose $x$ is adjacent to all
vertices of $\Gamma Y$.
Suppose $y$ is another vertex in $\Gamma X$.
Then, by Lemma~\ref{tools}$(ii)$, $y$ has fewer than $|\Gamma Y|-2$
neighbors in $\Gamma Y$, contradiction.
Similarly, if $|\Gamma Y|\geq 2$ each vertex in $\Gamma Y$ is adjacent to all but one vertices
of $\Gamma X$.
This implies that the subgraph induced by $\Gamma X \cup \Gamma Y$ is
$K_2$ or a complete bipartite graph with the edges of a perfect matching deleted.
So we find two possible block structures and quotient matrices for $A$:
\[
A=\left[
\begin{array}{cccc}
J-I_a   & J       & \1 & \0 \\
J       & J-I_b   & \0 & \1 \\
\1^\top & \0^\top & 0  & 1  \\
\0^\top & \1^\top & 1  & 0
\end{array}\right],\ 
Q=\left[\begin{array}{cccc}
a-1 & b   & 1 & 0 \\
a   & b-1 & 0 & 1 \\
a   & 0   & 0 & 1 \\
0   & b   & 1 & 0 \\
\end{array}\right],
\]
or
\[
A=\left[
\begin{array}{cccc}
J-I_a & J     & J     & O \\
J     & J-I_b & O     & J \\
J     & O     & O     & J-I_m \\
O     & J     & J-I_m & O
\end{array}\right],\ 
Q=\left[\begin{array}{cccc}
a-1 & b   & m   & 0 \\
a   & b-1 & 0   & m \\
a   & 0   & 0   & m-1 \\
0   & b   & m-1 & 0 \\
\end{array}\right],
\]
where $m = |\Gamma X|=|\Gamma Y|$. 
In the former case, $Q$ has eigenvalue 1 with multiplicity 1 for $(a,b)\in\{(4,4),(6,3)\}$, 
but none of the other 3 eigenvalues are equal to $\pm1$. 
In the latter case, $Q$ has eigenvalue 1 with multiplicity 1 for $(a,b,m)\in\{(3,3,8),(4,3,7),(4,4,6),(6,6,5),(8,5,5)\}$, 
but none of the other 3 eigenvalues are equal to $\pm1$.
For any other $a$, $b$, and $m$, neither quotient matrix has any eigenvalue equal to $\pm 1$ 
(this follows straightforwardly by solving $\det(Q+I)=0$ and $\det(Q-I)=0$).
Therefore the corresponding graphs are not in $\G$.
\qed
\\
{\bf Claim~2~} $a=b=2$.
\\[3pt]
{\bf Proof}.
First assume $a > b = 1$.
Then $\Gamma X$ contains no edges, because otherwise $C$ would not be maximal.
Consider $u\in Y$ and $x\in \Gamma X$.
Then $x$ is adjacent to all vertices in $\Gamma Y$ since otherwise interchanging $u$ and $x$
would give another maximal clique of size $c$ with fewer outgoing edges.
This implies that $u$ and $x$ have the same neighbors, contradiction.

Next assume $a > b = 2$.
We see that $a\geq 3$ implies that $\Gamma Y$ contains no edges, otherwise $G$ contains forbidden graph J.
Take a vertex $u\in X$ and a vertex $x\in\Gamma X$.
If $y$ and $z$ are distinct vertices in $\Gamma X$ both adjacent to $x$,
then the graph induced by $\{u,x,y,z\}\cup Y$ is a forbidden subgraph
(equal to J or containing D) of Lemma~\ref{tabu-l}.
Therefore $\Gamma X$ contains no intersecting edges.
Like before, forbidden graph~S implies that every vertex in $\Gamma X$ is adjacent to all, or all vertices but one in $\Gamma Y$.
Consider a vertex $x\in\Gamma X$ with no neighbors in $\Gamma X$ and a vertex $y\in Y$.
The neighborhood of $x$ is contained in that of $y$, but $d_y\leq d_x+2$, which contradicts Lemma~\ref{tools}$(ii)$.
We conclude that the graph induced by $\Gamma X$ is a disjoint union of edges.
Suppose $\{x,y\}$ is an edge in $\Gamma X$.
Then both $x$ and $y$ are adjacent to all vertices of $\Gamma Y$,
since otherwise interchanging $\{x,y\}$ with $Y$ would give another
clique in $G$ of size $c$ with fewer outgoing edges. Thus, every vertex of $\Gamma X$ is adjacent to every vertex of $\Gamma Y$.

Applying Lemma~\ref{tools}$(ii)$ to two vertices in $\Gamma Y$ 
yields a contradiction, thus $|\Gamma Y|=1$.
We find the following $A$ and $Q$:
\[
A=\left[
\begin{array}{ccc}
J-I_a   & J     & \0 \\
J       & R_{2m} & \1 \\
\0^\top & \1^\top & 0
\end{array}\right],\ 
Q=\left[\begin{array}{ccc}
a-1 & 2m & 0 \\
a   &     1 & 1 \\
0   & 2m & 0
\end{array}\right] .
\]
It follows straightforwardly that $Q$ has no eigenvalue equal to $\pm 1$.
This is a contradiction, and we conclude that $a=b=2$.
\qed
\\[3pt]
We have $a=b=2$.
By the same argument as above it follows that $\Gamma X$ only contains disjoint edges.
Forbidden graphs H and S imply that every vertex of $\Gamma X$ is adjacent to all, or all vertices but one of $\Gamma Y$.
Then, as before, a vertex in $\Gamma X$ with no neighbors in $\Gamma X$ and a vertex in $Y$ violate Lemma~\ref{tools}$(ii)$,
so $\Gamma X$ induces a disjoint union of edges.
Thus, every vertex of $\Gamma X$ must be adjacent to every vertex of $\Gamma Y$.
The same holds if $X$ and $Y$ are interchanged.
Thus we can conclude that $A$ is as follows:
\[ A=\left[\begin{array}{cc}
R_{2\ell} & J \\
J      & R_{2m}
\end{array}\right]
\]
with $\ell,m\geq 2$, where $2\ell=|\Gamma X|+2$ and $2m=|\Gamma Y|+2$. This leads to Case~$(iii)$.
\qed

\subsection{$\Omega$ is not empty}

Since $G$ is connected there exists an edge $\{x,z\}$ with $z\in\Omega$ and $x\in\Gamma X$, or $x\in\Gamma Y$. 
Assume $x\in\Gamma X$, take $u\in X$, and let $y$ be a neighbor of $z$ different from $x$. 
If $y\in\Gamma Y$, then the neighbor $v\in Y$ of $y$ together with $u$,
$x$, $y$, and $z$ induce a forbidden subgraph A or B from Lemma~\ref{tabu-l}.
Thus, $y\not\in\Gamma Y$ which means $y\in \Gamma X \cup \Omega$.

Assume that $|Y|\geq 2$. Let $v$ and $w$ be distinct vertices in $Y$.
If $y\in \Gamma X$, then $\{u,v,w,x,y,z\}$ induces a forbidden subgraph of type G or H. 
If $y\in \Omega$, then $\{u,v,w,x,y,z\}$ induces a forbidden subgraph of type K or M. Therefore $|Y|=1$.

Consider the set $Y'=Y\cup\Gamma X$ and let $Z$ be the set of vertices which are not in $X$ or $Y'$.
Then $|Y'|\geq 2$, since $Y$ and $\Gamma X$ are nonempty, and $Y'$ contains no edges, since otherwise $C$ wouldn't be maximal.
Therefore $X,Y'$ and $Z$ give the following block structure of $A$:
\[
A=\left[
\begin{array}{ccc}
J-I_a   & J     & O \\
J       & O & N \\
O & N^\top & M
\end{array}\right].
\]
Take three vertices $u\in X$, $x\in Y'$ and $y\in Y'$ with degrees $d_u$, $d_x$ and $d_y$, respectively.
Assume $d_x\leq d_y$, and consider the corresponding $3\times 3$ principal submatrix $S$ of $A^2-I$.
Then
\[
S=\left[
\begin{array}{ccc}
d_u-1 & a-1   & a-1 \\
a-1   & d_x-1 & d_{xy} \\
a-1   & d_{xy}& d_y-1
\end{array}\right] ,
\]
where $d_{xy}$ is the number of common neighbors of $x$ and $y$.
Write $S=(a-1)J + S'$, then
\[
S' = \left[\begin{array}{cc}d_u-a & \0^\top \\ \0 & T \end{array}\right] \mbox{ with }
T  = \left[\begin{array}{cc}d_x-a & d_{xy}-a+1 \\ d_{xy}-a+1 & d_y-a \end{array}\right].
\]
Note that $d_u>a$, and $d_x,d_y\geq a$. 
Forbidden graphs H and T imply that $y$ has at most two neighbors in $Z$ that are not neighbors of $x$. 
More precisely, if $y$ has two adjacent neighbors that are not neighbors of $x$, 
then these two neighbors of $y$ together with $x, y$ and two vertices in $X$ induce forbidden subgraph H. 
Otherwise, if $y$ has three neighbors that are not neighbors of $x$, 
then these three neighbors of $y$ form an independent set and together with $x, y$ and two vertices in $X$, induce forbidden subgraph T.
Thus, $d_x\leq d_y \leq d_{xy}+2\leq d_x+2$.

If $T$ is positive definite, then so are $S'$ and $S$, which contradicts rank$\, S \leq 2$.
Therefore $\det T = (d_x-a)(d_y-a)-(d_{xy}-a+1)^2 \leq 0$.
If $d_x=d_{xy}$ then $d_y\geq d_x+3$ (by Lemma~\ref{tools}$(ii)$), contradiction.
Also if $d_x\geq d_{xy}+1$, then $\det T > 0$ unless $d_x=d_y=d_{xy}+1$.
We conclude that $d_x=d_y=d_{xy}+1$, and we find the following two possible structures for $N$:
\[
N=[\, O\ \, J\!-\!I\ \, J\,]\ , \mbox{ or } N=[\, O\,\ I\,\ J\, ].
\]
Suppose a vertex $z\in Z$ has two neighbors $x$ and $y$ in $Z$.
Take three vertices $u$, $v$, and $w$ in $C$.
Then $\{u,v,w,x,y,z\}$ induce or contain a forbidden subgraph G, H, J, K, L, M, P, Q, C, or D from Lemma~\ref{tabu-l}.
So a vertex $x\in Z$ has at most one neighbor in $Z$, and since all vertices have degree at least two,
$z$ is adjacent to a vertex of $Y'$, hence $N=[\, J\!-\!I\ \, J\,]$ or $N=[\, I\,\ J\, ]$.
Partition $Z=Z_1\cup Z_2$ according to the structure on $N$, so that the vertices in $Z_2$ are adjacent to all
vertices of $Y'$.

Suppose $\{y,z\}$ is an edge in $Z$, and suppose there is a vertex $x\in Y'$ adjacent to $y$ but not to $z$. 
Take $u\in X$ and let $v\in Y'$ be a neighbor of $z$. 
Then $\{u,v,x,y,z\}$ induces a forbidden subgraph A. 
So $y$ and $z$ have the same set of neighbors in $Y'$, and hence $y,z\in Z_2$.
Take $w\in Z_1$. 
Then $w$ has no neighbor in $Z$, and hence has at least two neighbors in $Y'$ so $N=[\, J\!-\!I\ \, J\,]$. 
If $w\in Z_1$ and $z\in Z_2$, then every neighbor of $w$ is also a neighbor of $z$, but the degrees of $z$ and $w$ differ by two. 
This is impossible by Lemma~\ref{tools}$(ii)$.
Clearly $|Z_1|=|Y'|$ so $Z_1$ is not empty.
The conclusion is that $Z_2$ is empty and we find the following $A$ and $Q$:
\[
A=\left[
\begin{array}{ccc}
J-I_a & J     & O \\
J     & O     & J-I_m \\
O     & J-I_m & O
\end{array}\right]
\, , \
Q=\left[
\begin{array}{ccc}
a-1 & m     & 0 \\
a   & 0     & m-1 \\
0   & m-1 & 0
\end{array}\right]
,
\] where $m=|Y'|=|Z|$.
The matrix $Q$ has all three eigenvalues unequal to $\pm 1$, except when $(a,m)$ equals $(4,4)$ or $(3,5)$.
This leads to Case~$(vi)$.\qed

\section*{Acknowledgments} 
The work of the first and third author was partially supported by National Security Agency grant H98230-13-1-0267.

\end{document}